\documentclass[12pt]{article}
\usepackage{amsmath,amssymb}

\newcommand{\rest}[1]{_{|#1}}
\newcommand{\sD}{\mathcal D}
\newcommand{\Oh}{\mathcal O}
\newcommand{\PP}{\mathbb P}
\newcommand{\Ybar}{\overline{Y}}
\newcommand{\fie}{\varphi}

\newcommand{\Z}{\mathbb Z}

\DeclareMathOperator{\Gal}{Gal}
\DeclareMathOperator{\Qua}{Quad}
\newcommand{\broken}{\dasharrow}
\newcommand{\iso}{\cong}
\newtheorem{thm}{Theorem}[section]
\newtheorem{cor}[thm]{Corollary}
\newtheorem{prop}[thm]{Proposition}
\newtheorem{princ}[thm]{Principle}
\newtheorem{lem}[thm]{Lemma}

\newcommand{\step}[2]{\paragraph{\sc Step #1. \em#2}}

\newcommand{\pionealg}{\pi_1^{\mathrm{alg}}}
\newcommand{\F}{\mathbb F}
\DeclareMathOperator{\Sing}{Sing}
\newcommand{\Si}{\Sigma}

\newcommand{\1}{^{-1}}
\numberwithin{equation}{section}

\title{Campedelli surfaces with\\ fundamental group of order 8}
\author{Margarida Mendes Lopes\thanks{ The first author is a member of
the Centre for Mathematical Analysis, Geometry and Dynamical Systems,
Instituto Superior T{\'e}cnico, Lisboa. The second is a member of
G.N.S.A.G.A.--I.N.d.A.M. This research was partially supported by the
Italian project ``Geometria delle variet\`a algebriche e dei loro spazi
di moduli'' (PRIN COFIN 2006) and by FCT (Portugal) through program
POCTI/FEDER.} \and Rita Pardini${^*}$ \and Miles Reid }

\date{}

\begin{document}
\maketitle
\begin{abstract}
Let $S$ be a Campedelli surface (a minimal surface of general type with
$p_g=0$, $K^2=2$), and $\pi\colon Y\to S$ an etale cover of degree $8$.
We prove that the canonical model $\Ybar$ of $Y$ is a complete
intersection of four quadrics $\Ybar=Q_1\cap Q_2\cap Q_3\cap
Q_4\subset\PP^6$. As a consequence, $Y$ is the universal cover of $S$,
the covering group $G=\Gal(Y/S)$ is the topological fundamental group
$\pi_1S$ and $G$ cannot be the dihedral group $D_4$ of order $8$.
\newline {\em Mathematics Subject Classification (2000)}: 14J29.
\end{abstract}

\section{Introduction}

Let $Y$ be a minimal surface of general type with $K_Y^2=16$ and
$p_g=7$, $q=0$, having a free action by a group $G$ of order 8. Write
$\fie\colon Y\broken \Ybar\subset \PP^6$ for the $1$-canonical map, with
image $\Ybar$. We prove the following:

\begin{thm}\label{main}
The surface $\Ybar\subset\PP^6$ is the complete intersection of $4$ quadrics.
It is isomorphic to the canonical model of $Y$.
\end{thm}

Theorem~\ref{main} is known if $G=\Z_2^3$ by Miyaoka \cite{miyaoka},
Theorem~B; in this case there are four linearly independent diagonal
quadrics through $\Ybar$, which necessarily form a regular sequence
defining $\Ybar$. We thus assume throughout that $G$ is a group of order
8 and contains an element of order 4.

\begin{cor} Let $S$ be a Campedelli surface and $\pi\colon Y\to S$ an
etale cover of degree $8$. Then $Y$ is the universal cover of $S$ and
the covering group $G=\Gal(Y/S)$ is the topological fundamental group
$\pi_1S$. \end{cor}

\begin{cor}\label{dihedral} The dihedral group $D_4$ of order $8$ is not
the fundamental group of a Campedelli surface. \end{cor}

The proof of Theorem~\ref{main} consists of two parts, the first of
which is now quite standard (compare Reid \cite{Re2}, Naie \cite{Na},
Konno \cite{Ko}):

\begin{prop} \label{prop!1}
\begin{enumerate}
\renewcommand{\labelenumi}{(\roman{enumi})}
\item The canonical linear system $|K_Y|$ on $Y$ is free and defines a
morphism $\fie\colon Y\to\Ybar\subset \PP^6$ that is birational to its
image.
\item If $\Ybar$ is not a complete intersection of four quadrics, its
quadric hull (the intersection of all quadrics containing $\Ybar$) is a
$3$-fold $X$ of degree $4$, $5$ or $6$.
\item Moreover, in these three cases, $\Ybar$ is contained in a
hypersurface $F_d$ of $\PP^6$ not containing $X$, of degree $d=6,4,3$
respectively.
\end{enumerate}
\end{prop}

The second part analyses the possible cases $\Ybar\subset X$, with ad
hoc arguments involving the $G$-action to rule out each case; see
Section~\ref{s!allX}.

\subsection{The background}
A Campedelli surface is a surface $S$ of general type with $p_g=0$,
$K^2=2$. The algebraic fundamental group $\pionealg(S)$ classifies
finite etale covers $Y\to S$, and is the profinite completion of the
topological fundamental group $\pi_1S$. Results of Beauville \cite{Be}
and Reid \cite{Re1}, \cite{Re2} (see also Mendes Lopes and Pardini
\cite{3chi}) guarantee that $S$ has no irregular covers, and that an
etale cover $Y\to S$ has degree $\le9$. The reasons underlying
\cite{Re1}, Theorem~\ref{main} and all related results are as follows:

\begin{princ} \label{princ!16}
\renewcommand{\labelenumi}{(\arabic{enumi})}
\begin{enumerate}
\item The automorphism group $G$ acts on any intrinsically defined
feature of $Y$: for example, the base points or base $-2$-cycles of
$|K_Y|$ occur in multiples of $8$.

\item If a subgroup $H\subset G$ normalises a subscheme $Z\subset Y$,
its order $|H|$ divides the Euler characteristic $\chi(\Oh_Z)$; for
example, if $Y$ has an intrinsically defined genus $g$ pencil
$\psi\colon Y\to\PP^1$ and $H\subset G$ fixes $P\in\PP^1$ then $|H|$
divides $\chi(\Oh_F)=g-1$, where $F=\psi^*P$.
\end{enumerate}
\end{princ}

It seems most likely that all groups of order $\le9$ except the dihedral
groups of order 8 and 6 occur as $\pi_1S$. The case $|\pi_1S|=9$ was
treated in detail in Mendes Lopes and Pardini \cite{MP3}. Here we treat
$|\pi_1S|=8$, patching up the incomplete manu\-script \cite{Re2}. Naie
\cite{Na} obtained similar results for $|\pi_1S|=6$ using similar
methods. Campedelli surfaces with $\pi_1=\Z/8$ and $\Z/2\oplus\Z/4$ are
contained in passing in Barlow \cite{Ba}. Beauville \cite{Be2}
constructs a family of Calabi--Yau 3-folds with $\pi_1$ the quaternion
group $H_8$, and Campedelli surfaces with the same $\pi_1$ are obtained
by taking the unique invariant section $X_1=0$ of this.

\subsection{Representations of $G$ and proof of
Corollary~\ref{dihedral}} Let $Y\to S$ be the universal cover of a
Campedelli surface with group $G$. Then $G$ acts naturally on $H^0(K_Y)$
and $H^0(2K_Y)$. Since the $G$-action is free,  $H^0(K_Y)$ is the
regular representation of $G$ minus the trivial rank~1 representation,
and $H^0(2K_Y)$ is three times the regular representation (for example,
by \cite[Corollary~8.6]{ypg}). Finally, the $G$-equivariant
multiplication map
\begin{equation}\label{eq!mult}
S^2H^0(K_Y)\to H^0(2K_Y)
\end{equation}
is surjective by Theorem~\ref{main}.

These remarks allow one to show that the group $G$ is not the dihedral
group, and to describe explicitly $Y$ and the $G$-action for all the
remaining groups of order 8.

Let $G=D_4$ be the dihedral group of order 8. Write $1$ for the trivial
rank~1 representation, and $\rho$ for the sole irreducible rank~2
representation; let $\chi_1:=\bigwedge^2\rho$, $\chi_2$ and $\chi_3$ be
the remaining rank~1 representations. Then
\begin{equation}
\begin{aligned}
H^0(K_Y)&=\chi_1\oplus \chi_2\oplus \chi_3\oplus \rho^{\oplus 2}, \\[4pt]
H^0(2K_Y)&=1^{\oplus 3}\oplus \chi_1^{\oplus 3}\oplus \chi_2^{\oplus
3}\oplus \chi_3^{\oplus 3} \oplus \rho^{\oplus 6}.
\end{aligned}
\end{equation}
Using the decomposition of $H^0(K_Y)$, one computes:
\begin{equation}
S^2H^0(K_Y)=1^{\oplus 6}\oplus \chi_1^{\oplus 2}\oplus \chi_2^{\oplus
4}\oplus \chi_3^{\oplus 4}\oplus \rho^{\oplus 6}.
\end{equation}
Clearly the equivariant map \eqref{eq!mult} cannot be surjective. This
contradicts Theorem~\ref{main} and proves Corollary~\ref{dihedral}.

\section{Proof of Proposition~\ref{prop!1}}

The canonical map $\fie\colon Y\to\PP^6$ is a morphism by Ciliberto,
Mendes Lopes and Pardini \cite[Proposition~5.2]{cmp} and is birational
to its image $\Ybar$ by \cite[Proposition~5.3]{cmp}. Thus $\Ybar$ is an
irreducible surface of degree~16. Since
\[
\dim S^2H^0(Y,K_Y)=\binom82=28 \quad\hbox{and}\quad
h^0(Y,2K_Y)=\chi(\Oh_Y)+K^2_Y=24,
\]
the multi\-plication map
$S^2H^0(Y,K_Y)\to H^0(Y, 2K_Y)$ has kernel of dimension $\ge4$; that is,
$\Ybar$ is contained in at least 4 linearly independent quadrics.

Let $Q_1,Q_2,Q_3,Q_4$ be four linearly independent quadrics through
$\Ybar$. We are home if $\Ybar$ is an irreducible component of
$Q_1\cap\cdots\cap Q_4$. For in turn, if any of $Q_1$ or $Q_1\cap Q_2$
or $Q_1\cap Q_2\cap Q_3$ or $Q_1\cap\cdots\cap Q_4$ is reducible, then
$\deg\Ybar<16$. This is impossible, so $\Ybar=Q_1\cap\cdots\cap Q_4$ is
a complete intersection of 4 quadrics. Then $\Ybar$ is Gorenstein with
$K_{\Ybar}=\Oh_{\Ybar}(1)$ and $K_Y=\fie^*K_{\Ybar}$. Therefore it has
canonical singularities and is the canonical model of $Y$.

Write $\Qua(\Ybar)\subset\PP^6$ for the {\em quadric hull} of $\Ybar$,
the intersection of all the quadrics through $\Ybar$, following
\cite{Re3} and Konno \cite{Ko}. The alternative to $\Ybar$ a complete
intersection of four quadrics is that $\Qua(\Ybar)$ has a component $X$
strictly containing $\Ybar$. Then $X$ is a $3$-fold of degree $4$, $5$
or $6$ and is the unique component of $\Qua(\Ybar)$ containing $\Ybar$.

Indeed, by elementary inequalities due to Castelnuovo, an irreducible
$m$-fold $X$ spanning $\PP^N$ is contained in at most
\begin{equation}
\binom{N-m+2}2-\min\{\deg X,2(N-m)+1\}
\label{eq!ineq}
\end{equation}
linearly independent quadrics. See for example the discussion in
\cite{Re3} or \cite[Corollary~1.5]{Ko}. The equality
$X=\Qua(\Ybar)$ follows by \cite[Corollary~2.6]{Ko}. The estimate on
$d$ follows from (\ref{eq!ineq}) or by \cite[Proposition~1.3]{Ko}.

Finally, in the three cases for $d$, crude estimates give that the
restriction map
\[
H^0(\PP^6,\Oh(k)) \to H^0(\Oh_X(k))
\]
has rank
\begin{equation}
\begin{tabular}{lrl}
$=252$ & for $d=4$, $k=6$, & whereas $h^0(6K_Y)=248$;\\[4pt]
$\ge105$ & for $d=5$, $k=4$, & whereas $h^0(4K_Y)=104$;\\[4pt]
$\ge58$ & for $d=6$, $k=3$, & whereas $h^0(3K_Y)=56$
\end{tabular}
\end{equation}
(compare \cite{Re3} and \cite[Lemma~1.8]{Ko}). This proves
Proposition~\ref{prop!1}.

\section{Proof of Theorem~\ref{main}}\label{s!allX}

We exclude the cases of Proposition~\ref{prop!1}, (ii) by
studying the $G$-action on $\Ybar\subset X$, treating separately the cases $\deg
X=4$, $5$ or $6$. In any case, $X$ is linearly normal, since
$Y\to\Ybar\subset\PP^6$ is given by the complete canonical system, and is
regular, since $Y$ is.

\subsection{$G$-invariant linear systems on $Y$}

The following lemmas group together a number of restrictions on
$G$-invariant linear systems on $Y$, that we use several times in what
follows. Their proofs are applications of Principle~\ref{princ!16}.

\begin{lem}\label{involution}
A $G$-invariant linear system $|D|$ on $Y$ with $D^2=2$ has a fixed
part.
\end{lem}

\begin{pf} Assume by contradiction that $|D|$ has no fixed part. Since
$G$ acts on the base locus of $|D|$, $D^2=2$ implies $|D|$ is free.
Hence $|D|$ defines a $G$-equivariant $2$-to-$1$ morphism $Y\to\PP^2$.
Since we assume that $G$ has an element of order 4, this contradicts
Beauville \cite[Corollary 5.8]{Be}. \end{pf}

\begin{lem}\label{pencil} Let $|F|$ be a $G$-invariant pencil on $Y$
with $K_YF\le 8$. Then $|F|$ is free and $K_YF=8$. \end{lem}

\begin{pf} Since $K_YF\le 8$, the index theorem gives $F^2\le
(K_YF)^2/K_Y^2\le 4$. Now $F^2$, equal to the degree of the base locus
of $|F|$, is divisible by 8 by Principle~\ref{princ!16}, so $F^2=0$ and
$|F|$ is free. If $K_YF<8$, the general $F\in|F|$ is nonsingular of
genus $g\le 4$, contradicting \cite[Lemma~2.2]{cmp}, so $K_YF=8$.
\end{pf}

\begin{prop}\label{system}
Let $\sD\subset |K_Y|$ be a $G$-invariant subsystem of projective
dimension $\ge3$. Then one of the following holds:
\begin{enumerate}
\renewcommand{\labelenumi}{(\arabic{enumi})}
\item $\sD$ is free; or
\item $\sD$ has base locus consisting of $8$ transversal base points.
\end{enumerate}
In particular, $\sD$ is not composed with a pencil.
\end{prop}

\begin{pf} If $\sD$ has a nonzero fixed part $Z$, write $\sD=M+Z$ with
mobile $|M|$. The $G$-action takes $Z$ to itself, so $Z$ is the pullback
from $S=Y/G$ of a divisor $Z_0$ that satisfies $K_SZ_0\equiv Z_0^2$
mod~2; therefore $MZ=(K_Y-Z)Z$ is divisible by 16. Connectedness of
canonical divisors gives $MZ>0$, and thus $16\le MZ\le K_YM\le
K^2_Y=16$. We get:
\[
M^2=K_YZ=0\quad\hbox{and}\quad K_YM=16, \quad Z^2=-16.
\]

Since $|M|$ is mobile and $M^2=0$, it is contained in a multiple of a $G$-invariant
free pencil, say $|M|\subset |nF|$ with $K_YF=16/n$; Lemma~\ref{pencil} implies
$n\le2$. But $n=\dim|nF|\ge\dim\sD\ge3$, a contradiction.

Therefore $\sD$ has no fixed part. Since $G$ acts on the base scheme of
$\sD$, the number $\nu$ of base points is divisible by 8. If $\nu>8$ or
$\nu=8$ and the base points are not transversal, two curves of $\sD$
have no free intersections, hence $\sD$ is composed with a pencil. Write
$D=nF$, with $|F|$ a $G$-invariant pencil and $n\ge3=\dim\sD$. Then
$F^2=0$ by Lemma~\ref{pencil}, contradicting $16=D^2=n^2F^2$.
\end{pf}

\subsection{The case $\deg X=4$}
In this case, by Fujita \cite{Fuj1}, $X$ is either a quartic scroll
$\F(a,b,c)$ with $a+b+c=4$, or the cone over the Veronese surface
$V_4\subset \PP^5$. By Proposition~\ref{prop!1}, (iii), there is a
sextic hypersurface containing $\Ybar$ and not containing $X$.

If $X$ is a scroll, the birational transform of its unique ruling by
planes is a $G$-invariant pencil $|F|$ on $Y$ with $K_YF\le 6$,
contradicting Lemma~\ref{pencil}.

If $X$ is the cone over $V_4$, the linear subsystem $\sD\subset |K_Y|$
formed by hyperplanes through its vertex define a $G$-equivariant map
$\psi\colon Y\to V_4$. By Proposition~\ref{system}, $\sD$ is either free
or has 8 simple base points. In the latter case, $\deg\psi=2$
contradicts \cite[Corollary 5.8]{Be}, as in Lemma~\ref{involution}. So
$\sD$ is free, and $\psi\colon Y\to V_4\iso\PP^2$ is a morphism of
degree~4. The $G$-action on $\PP^2$ fixes some point $P\in\PP^2$ by
Lemma~\ref{l!fPP2} below, whereas $\psi\1P$ consists of $\le4$ points or
trees of $-2$-curves, on which $G$ cannot act freely. This is a
contradiction.

\begin{lem}\label{l!fPP2} Let $G$ be a group of order $2^r$ acting on
$\PP^2$. Then there is a point $P\in\PP^2$ fixed by the whole of $G$.
\end{lem}
\begin{pf}
Indeed, $G$ has nontrivial centre, so a central element $g$ of order~2.
The action of $g$ on $\PP^2$ must fix an isolated point $P$ and a line
$L$. For any $h\in G$, by the conjugacy principle, the element $hgh\1$
is an involution with isolated fixed point $h(P)$. But $hgh\1=g$, so
that $h(P)=P$.
\end{pf}

\subsection{The case $\deg X=5$}

This is the hard case of the proof, and we break it into several steps.

\step1{$X$ is a normal del Pezzo variety with $K_X=\Oh_{X}(-2)$.}
Recall from the start of the proof that we assume that $X$ is linearly
normal and regular. By \cite[Theorem~2.1]{Fuj3} (or \cite{Fuj2} in the
nonsingular case) $X$ is either a normal del Pezzo variety of index~2 or
a cone from a point vertex over a (weak) del Pezzo surface $V_5\subset
\PP^5$. If $X$ is a cone, the subsystem $\sD\subset |K_Y|$ given by
hyperplanes through its vertex defines a $G$-equivariant map $\psi\colon
Y\to V_5$. By Proposition~\ref{system}, $\psi$ is onto the surface
$V_5$, and
\[
\deg V_5\cdot\deg\psi=5\deg \psi=8 \hbox{ or } 16
\]
provides a contradiction.

\step2{$\Ybar\cap\Sing X$ is a finite set.}
If $\Sing X$ is positive dimensional, it contains a single line $L$
(\cite[Theorem~2.7]{Fuj3}). Apply Proposition~\ref{system} to the
subsystem $\sD\subset|K_Y|$ given by hyperplanes of $\PP^6$ through $L$;
then $\sD$ has no fixed part, so $L$ is not contained in $\Ybar$.

\step3{The general section $C$ of $\Ybar$ is nonsingular.}
Let $\Si$ be a general hyperplane section of $X$ and set
$C=\Si\cap\Ybar$. The surface $\Si$ is a (possibly singular) del Pezzo
surface of degree 5, nonsingular along $C$ by Step~2, so that $C$ is a
Cartier divisor on $\Si$. Write $A=-K_\Si=\Oh_\Si(1)$ for the
restriction of a hyperplane to $\Si$. Since $AC=-K_{\Si}C=16$, the index
theorem gives $C^2\le(AC)^2/A^2=\frac{256}5$, so $C^2\le 51$. The curve $C$ is the
birational image of a general canonical curve of $Y$, so has geometric
genus 17. On the other hand, the arithmetic genus of $C\subset\Si$ is
given by $2p_aC-2=C^2+K_\Si C=C^2-16$. There are thus two possibilities:
\begin{enumerate}
\renewcommand{\labelenumi}{(\alph{enumi})}
\item $C^2=48$ and $C$ is nonsingular, or
\item $C^2=50$ and $C$ has a single node or cusp.
\end{enumerate}

If case (b) holds for the general hyperplane section of $\Ybar$, the
codimension~1 part of the singular locus of $\Ybar$ is a line $L$,
necessarily invariant under the action of $G$. The system of hyperplanes
through $L$ then give the same contradiction to Proposition~\ref{system}
as in Step~2.

\step4{Conclusion of the proof.} We continue to use the notation of Step~3. The
canonical class of $C$ calculated on $Y$ is
$K_C=(K_Y+C)\rest{C}=\Oh_C(2A)$. Calculated on $\Si$, it is
$(K_\Si+C)\rest{C}=\Oh_C(-A+C)$. Therefore the Cartier divisor $D=C-3A$
on $\Si$ restricted to $C$ is linearly equivalent to zero. Consider the
exact sequence of sheaves on $\Si$:
\[
0\to\Oh_{\Si}(-3A)\to \Oh_{\Si}(D)\to\Oh_C\to 0.
\]
Since $H^1(\Si,-3A)=0$ by Kodaira vanishing, or by well known results on
del Pezzo surfaces, it follows that $h^0(\Oh_\Si(D))=1$, so $D$ is a
Cartier divisor linearly equivalent to an effective divisor.

Now $-K_\Si D=AD=1$, and $D^2=48-96+45=-3$. This is a contradiction.
Indeed, $AD=1$ and $A$ very ample implies that $D$ is a line on $\Si$.
But then $D$ is nonsingular, and because it is a Cartier divisor, $\Si$
is nonsingular near $D$, so $D^2=-1$.

\subsection{The case $\deg X=6$}

Assume $\deg X=6$. By Proposition~\ref{prop!1} and its proof, the linear
system of cubics of $\PP^6$ containing $\Ybar$ restricts on $X$ to a
positive dimensional linear system $|N|$ of surfaces of degree~2. Now
$X$ is not ruled by planes (because it is linearly normal of degree 6
and regular), so that the moving part of $|N|$ must be a pencil of
quadrics.

The birational transform of $|N|$ on $Y$ is then a $G$ invariant pencil
$|F|$ with $K_YF\le 6$ and contradicts Lemma~\ref{pencil}.

\bigskip

\bigskip
\noindent {\em Authors' addresses}:

\medskip

\noindent Margarida Mendes Lopes,\\
Departamento de Matem\'atica,\\ Instituto Superior T\'ecni\-co,
Universidade T\'ecnica de Lisboa,\\
Av.~Rovisco Pais, 1049-001 Lisboa, Portugal

\noindent {\it e-mail}: mmlopes@math.ist.utl.pt

\bigskip

\noindent Rita Pardini,\\Dipartimento di Matematica,\\
Universit\`a di Pisa,\\
Largo B. Pontecorvo, 5,
56127 Pisa, Italy

\noindent {\it e-mail}:
pardini@dm.unipi.it

\bigskip

\noindent Miles Reid\\
Math Institute,\\
Univ. of Warwick, \\
Coventry CV4 7AL

\noindent {\it e-mail}:
Miles.Reid@warwick.ac.uk
\end{document}